\documentclass[12pt]{article}
\usepackage{graphicx,psfrag,epsfig}
\usepackage{graphicx,psfrag,epsfig, color}
\usepackage{amssymb,amsmath,amscd,amsthm}
\usepackage{graphicx,psfrag,epsfig}

\usepackage{graphicx}
\usepackage[active]{srcltx}

\newtheorem{theorem}{Theorem}[section]

\newtheorem{lemma}[theorem]{Lemma}

\date{}

\begin{document}

\date{}
\title{On stochastic perturbations of slowly changing dynamical systems}
\author{ M. Freidlin\footnote{Dept of Mathematics, University of Maryland,
College Park, MD 20742, mif@math.umd.edu}, L.
Koralov\footnote{Dept of Mathematics, University of Maryland,
College Park, MD 20742, koralov@math.umd.edu}
} \maketitle

\begin{abstract}
In this paper we consider a diffusion process obtained as a small random perturbation of a dynamical system attracted to a stable equilibrium point. The drift and the diffusive perturbation are assumed to evolve slowly in time. We describe the asymptotics of the
time it takes the process to exit a given domain and the limiting distribution of the exit point.
\end{abstract}



\section{Introduction and the formulation of the main result}


Suppose that the state of a system is determined by a $d$-dimensional vector. Assume that without any perturbations the system is
situated at a point $O \in \mathbb{R}^d$, and that it is pushed back to $O$ when perturbed by noise, so that the evolution of the perturbed system is described by the equation
\[
d X^{x,\varepsilon}_t = b(  X^{x,\varepsilon}_t) d t+\varepsilon
\sigma ( X^{x,\varepsilon}_t) d W_t,~~X^{x,\varepsilon}_1 = x \in
\mathbb{R}^d.
\]
Here $b$ is a smooth vector field with an asymptotically stable equilibrium $O$ and $\sigma$ is a diffusion matrix. The vector field provides a repairing mechanism that returns the system close to the equilibrium.
Observe that the initial position of the process is prescribed at time $t =1$ rather than the usual $t=0$ for reasons that will become clear later.

Consider  a domain $D$ such that $O \in D$ and all the points of $\overline{D}= D \cup \partial D$ are attracted to $O$ for the unperturbed system (with $\sigma \equiv 0)$. Let $\tau^{x, \varepsilon}$ be the first time when $X^{x, \varepsilon}_t$ reaches $\partial D$. This can be viewed as the life span of our system - the system `dies' when $X^{x, \varepsilon}_t$ exits $D$. The asymptotics of the life span and the location of the point where $X^{x,\varepsilon}_t$ exits $D$ can be explicitly expressed through the action functional corresponding to the diffusion process (see \cite{FW}). In particular, $\lim_{\varepsilon \downarrow 0} (\varepsilon^2 \ln \mathrm{E}\tau^{x, \varepsilon}) = V > 0$, where $V$ depends on the domain $D$ and on the coefficients. The system also has the following renewal property: for every $x \in D$ and  $t(\varepsilon)$, the distribution of $\tau^{x, \varepsilon} - t(\varepsilon)$ conditioned on survival till time $t(\varepsilon)$ is asymptotically equivalent to the distribution of $\tau^{x, \varepsilon}$ (both are approximately exponential distributions with parameter $\exp(V/\varepsilon^2)$). This property can be interpreted as the lack of ageing in the system.

In certain applications, it is natural to allow the system to age. For instance, the repairing capabilities may degrade with time. Since in the absence of ageing the life span of the system is of order $\exp({\rm const}/\varepsilon^2)$, it is natural to assume that the degradation occurs slowly, i.e., also at exponential time scales. We allow for two ageing mechanisms - one due to the slow evolution of coefficients in time and the other due to the dependence on a slowly changing random process. This leads us to consider the following stochastic process.

Let $\xi_\lambda$, $\lambda \geq 0$,  be a continuous time Markov chain on the state space $S = \{1,...,s\}$. We choose the right-continuous modification of $\xi_\lambda$.
Consider the diffusion process
\[
d X^{x,\varepsilon}_t = b(  X^{x,\varepsilon}_t, \varepsilon^2 \ln t, \xi_{\varepsilon^2 \ln t}) d t+\varepsilon
\sigma ( X^{x,\varepsilon}_t,\varepsilon^2 \ln t, \xi_{\varepsilon^2 \ln t}) d W_t,~~X^{x,\varepsilon}_1 = x \in
\mathbb{R}^d.
\]
Here $\varepsilon > 0$ is a small parameter, $W_t$ is a Wiener
process in $ \mathbb{R}^d$, independent of $\xi_\lambda$. The coefficients $b(\cdot, \cdot, k)$, $1 \leq k \leq s$, and $\sigma(\cdot, \cdot, k)$
are assumed to be bounded, continuous, and Lipschitz continuous in the spatial variable (with the Lipschitz constant that doesn't depend on $\lambda$). The diffusion matrix $a(x, \lambda, k)
= (a_{ij}(x, \lambda, k)) = \sigma(x, \lambda, k) \sigma^*(x, \lambda, k)$ is assumed to be uniformly
positive definite. The initial position of the process is prescribed at time $t =1$ rather than the usual $t=0$ in order to avoid the large negative values of the logarithm inside the coefficients.

Let $D \subset \mathbb{R}^d$ be a bounded domain with smooth boundary. We assume that there is a point $O \in D$ (the equilibrium) and
$r ,c > 0$ such that
\[
(b(x, \lambda, k), x-O) \leq -c|x-O|^2
\]
whenever $x$ is in the $r$-neighborhood of $O$, $\lambda \geq 0$, and $k \in S$, and that
\[
(b(x, \lambda, k), n(x)) \leq -c,~~~x \in \partial D,
\]
where $n(x)$ is the outward unit normal vector.
Moreover, we assume that for each $\lambda$
every solution of the equation $x'(t) = b(x(t), \lambda, k)$ that starts in $\overline{D}$ enters the $r$-neighborhood of $O$ in time that is shorter than $c^{-1}$.

Recall that $\tau^{x, \varepsilon}$ is the first time when $X^{x, \varepsilon}_t$ reaches $\partial D$. We'll be interested in the asymptotics of $\tau^{x, \varepsilon}$ and the limiting distribution of $X^{x, \varepsilon}_{\tau^{x, \varepsilon}}$. First, consider the following family of auxiliary problems. For each $\lambda \geq 0$ and $k \in S$, let  $Y^{x,\varepsilon,\lambda,k}_t$ be the solution of
\[
d Y^{x,\varepsilon,\lambda,k}_t = b(  Y^{x,\varepsilon,\lambda,k}_t, \lambda, k) d t+\varepsilon
\sigma (  Y^{x,\varepsilon,\lambda,k}_t, \lambda, k) d W_t,~~Y^{x,\varepsilon,\lambda,k}_0 = x \in
\mathbb{R}^d.
\]
This is a diffusion process with time-independent coefficients, and the usual action functional can be defined:
\[
S_{0, T}^{\lambda,k} (\varphi) = \frac{1}{2}\int_{0}^{T} \sum_{i,j = 1}^d
a^{ij}(\varphi_t,\lambda,k)(\dot{\varphi}^i_t -
b_i(\varphi_t,\lambda,k))(\dot{\varphi}^j_t - b_j(\varphi_t,\lambda,k)) d t,~~T \geq
0,
\]
for absolutely continuous $\varphi \in C([0,T], \mathbb{R}^d)$, while $ S_{0, T}^{\lambda,k}(\varphi) =
+\infty$ for $\varphi$ that are not absolutely continuous. Here
$a^{ij}$ are the elements of the inverse matrix, that is $a^{ij} =
(a^{-1})_{ij}$.  The quasi-potential is defined as
\[
V^{\lambda,k}(x) =  \inf_{T, \varphi} \{ S_{0,T}^{\lambda,k}(\varphi): \varphi \in
C([0,T], \mathbb{R}^d), \varphi(0) = O, \varphi(T) = x \},~~x
\in \mathbb{R}^d.
\]
Let
\[
M^{\lambda, k} = \inf_{x \in \partial D}  V^{\lambda,k}(x) .
\]
Obviously, this is a continuous function of $\lambda$ with values in $(0,\infty)$.
\\
\\
{\bf Assumption 1.} We assume that the equation $M^{\lambda, k} = \lambda$ has a unique solution for each $k \in S$.
\\

The solution will be denoted by $m^k$. The case when the equation has finitely many solutions could also be considered without major additional
difficulties, but would require more complicated notations.
\\
\\
{\bf Assumption 2.} We assume, for brevity,  that for each $k \in S$ the infimum in $ \inf_{x \in \partial D}  V^{m^k,k}(x)$ is achieved in a single point of the boundary.
\\

\begin{figure}
  \begin{center}
     \includegraphics[height=4.2in, width= 4.9in,angle=0]{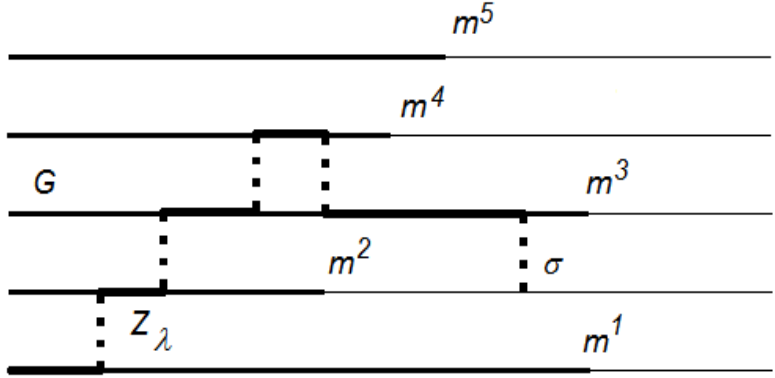}
  \end{center}
    \vskip -130pt
  \caption{A trajectory of the process $Z_\lambda$.}
    \label{stsp}
\end{figure}

The point for which the infimum is achieved will be denoted by $x^k$. Let $I_k = \{k\} \times [0, m^k)$, and $G = \bigcup_{k =1}^s I_k$. Thus $G$ can
be viewed as a union of $k$ disjoint segments of lengths $m^1,...,m^s$. It is a subset of the larger space $\overline{G} = S \times [0,\infty)$  (see Figure \ref{stsp}). Define the following Markov process on $\overline{G}$:
\[
Z_\lambda = (\xi_\lambda, \lambda),~~\lambda \geq 0,
\]
and let $\sigma = \inf\{\lambda: Z_\lambda \notin G\}$ be the first time when the process leaves $G$.

Now we can formulate the main result.
\begin{theorem} \label{mtee} Under the above assumptions, for each $x \in D$ the distribution of $\varepsilon^2 \ln \tau^{x, \varepsilon}$ converges, as $\varepsilon \downarrow 0$, to the distribution of $\sigma$. The distribution of $X^{x,\varepsilon}_{\tau^{x, \varepsilon}}$ converges to the distribution of $x^{\xi_\sigma}$.
\end{theorem}

Fix a number $\Lambda > \max (m^1,...,m^s)$ and a right-continuous function $z_\lambda$, $\lambda \in [0,\Lambda]$, taking values in $S$, with a finite number of jumps.  Together with $z_\lambda$, we can consider the function $(z_\lambda,\lambda)$ with values in $\overline{G} = S \times [0,\infty)$. This can be viewed as a trajectory of the process $Z_\lambda$. Given such a function $z$, we can define the process
\begin{equation} \label{perturb2}
d X^{x,\varepsilon,z}_t = b(  X^{x,\varepsilon,z}_t, \varepsilon^2 \ln t, z_{\varepsilon^2 \ln t}) d t+\varepsilon
\sigma (  X^{x,\varepsilon,z}_t, \varepsilon^2 \ln t, z_{\varepsilon^2 \ln t}) d W_t,~~X^{x,\varepsilon,z}_{1} = x \in
\mathbb{R}^d.
\end{equation}
This is different from the process $X^{x,\varepsilon}$ in that now the trajectory of the underlying Markov process is considered fixed.
Let $\tau^{x, \varepsilon,z}$ be the first time when $X^{x, \varepsilon,z}_t$ reaches $\partial D$.
Let $\sigma^z$ be the first time when $(z_\lambda,\lambda)$ leaves $G$, which corresponds to the stopping time $\sigma$ for the process $Z$. Observe that the probability of $\xi_\lambda$ experiencing a jump at any of the points $m_1,...,m_s$ is zero and that the number of jumps on $[0,\Lambda]$ is finite with probability one.
Therefore, by conditioning on a trajectory of the process $\xi_\lambda$, we immediately reduce Theorem~\ref{mtee} to the following.
\begin{lemma} \label{mnne}
Assume that the function $z_\lambda$ is right-continuous with finitely many jumps at points $\lambda_1,...,\lambda_k$ and that none of the jumps happens at any of the points $m_1,...,m_s$. Then for each $x \in D$ the distribution of $\varepsilon^2 \ln \tau^{x, \varepsilon,z}$ converges, as $\varepsilon \downarrow 0$, to $\sigma^z$. The distribution of $X^{x,\varepsilon,z}_{\tau^{x, \varepsilon,z}}$ converges to $x^{z_{\sigma^z}}$.
\end{lemma}

In the next section we gather some facts about processes with time-dependent coefficients that will be needed for the proof of Lemma~\ref{mnne}. The lemma itself will be proved in Section~\ref{prml}. Finally, in Section~\ref{exa} we briefly discuss a situation in which the vector field $b$ is equal to zero near  $ \partial D$.

\section{On diffusion processes with time-dependent coefficients}

Before we proceed with the proof of Lemma~\ref{mnne}, let us briefly discuss diffusion processes
whose coefficients are time-dependent, but are close to
functions that do not depend on time. For $T > 0$ and $\varphi,
\psi \in C([0,T], \mathbb{R}^d)$, we define $\rho_T(\varphi, \psi)
= \sup_{t \in [0,T]}|\varphi(t) - \psi(t)|$.

Let $\widetilde{a}^\varepsilon(t,x)$, $\overline{a}(x)$,  be uniformly positive
definite symmetric $d \times d$ matrices whose elements
$\widetilde{a}_{ij}^\varepsilon$, $\overline{a}_{ij}$, are continuous in $(t,x)$ ($x$ in case of $\overline{a}$)
and  Lipschitz continuous in $x$ with a Lipschitz constant $L$.  Assume that there are positive $k$ and $K$ such that
\[
k |\xi|^2 \leq \sum_{i,j =1}^d \widetilde{a}^\varepsilon_{ij}(t,x) \xi_i \xi_j \leq K |\xi|^2,~~~t \geq 0, x \in \mathbb{R}^d,  \xi \in \mathbb{R}^d.
\]
\[
k |\xi|^2 \leq  \sum_{i,j =1}^d \overline{a}^\varepsilon_{ij}(x) \xi_i \xi_j \leq K |\xi|^2,~~~ x \in \mathbb{R}^d,  \xi \in \mathbb{R}^d.
\]
Let $\widetilde{\sigma}^\varepsilon$, $\overline{\sigma}$ be a square matrices such that
$\widetilde{a}^\varepsilon = \widetilde{\sigma}^\varepsilon
(\widetilde{\sigma}^\varepsilon)^*$, $\overline{a} = \overline{\sigma}
(\overline{\sigma})^*$. We choose
$\widetilde{\sigma}^\varepsilon$, $\overline{\sigma}$ in such a way that
$\widetilde{\sigma}_{ij}^\varepsilon$ are continuous in
$(t,x)$, while $\widetilde{\sigma}_{ij}^\varepsilon$ and $\overline{\sigma}_{ij}$ are Lipschitz continuous in~$x$. Let $\widetilde{b}^\varepsilon(t,x)$, $\overline{b}(x)$, be continuous vector fields, Lipschitz continuous in $x$ with Lipschitz constant $L$.

Let $\widetilde{X}^{x,\varepsilon}_t$ and $\overline{X}^{x,\varepsilon}_t$ satisfy $\widetilde{X}^{x,
\varepsilon}_0 = \overline{X}^{x,\varepsilon}_0 =  x$ and
\[
d \widetilde{X}^{x,\varepsilon}_t = \widetilde{b}(t,
\widetilde{X}^{x,\varepsilon}_t) d t +\varepsilon
\widetilde{\sigma}^\varepsilon (t,
\widetilde{X}^{x,\varepsilon}_t) d W_t,~~d \overline{X}^{x,\varepsilon}_t = \overline{b}(\overline{X}^{x,\varepsilon}_t) d t +\varepsilon
\overline{\sigma} (\overline{X}^{x,\varepsilon}_t) d W_t.
\]
We will assume that the coefficients of the process
$\widetilde{X}^{x,\varepsilon}_t$ are close to those of
$\overline{X}^{x,\varepsilon}_t$. Namely, let us assume that
\[
 \sup_{(t,x) \in \mathbb{R}^+
\times \mathbb{R}^d} | \widetilde{b}_{i}^\varepsilon(t, x)-
\overline{b}_{i}(x)| \leq \varkappa,~~~
 \sup_{(t,x) \in \mathbb{R}^+
\times \mathbb{R}^d}   | \widetilde{a}_{ij}^\varepsilon(t, x)-
\overline{a}_{ij}(x)| \leq \varkappa,
\]
where $\varkappa$ is small.

The reason to introduce these processes is that we would like to study
the behavior of the process $ X^{x,\varepsilon,z}_t$ given by
(\ref{perturb2}) on such time intervals where the
coefficients do not change much with time. Thus it is convenient to consider
a generic process whose coefficients are close to
functions that don't depend on time.  Let $\overline{S}_{0,T}$ be the action functional for the processes $\overline{X}^{x, \varepsilon}_t$, that is
\[
\overline{S}_{0, T} (\varphi) = \frac{1}{2}\int_{0}^{T} \sum_{i,j = 1}^d
\overline{a}^{ij}(\varphi_t)(\dot{\varphi}^i_t -
\overline{b}_i(\varphi_t))(\dot{\varphi}^j_t - \overline{b}_j(\varphi_t)) d t,~~T \geq
0,
\]
for absolutely continuous $\varphi \in C([0,T], \mathbb{R}^d)$, while $ \overline{S}_{0, T}(\varphi) =
+\infty$ for $\varphi$ that are not absolutely continuous.
The next two lemmas show that $\overline{S}$ serves a purpose similar
to the action functional for the processes
$\widetilde{X}^{x,\varepsilon}_t$, even though the diffusion
coefficients for the process are time-dependent.
\begin{lemma} \label{ler1}
Suppose that
$\widetilde{a}^\varepsilon$, $\overline{a}$, ${\widetilde{b}}^\varepsilon$, and $\overline{b}$ are as above, and positive
constants  $k$, $K$, and $L$ are fixed. For every $\delta$, $\gamma$
and $C$ there exist $\varkappa > 0$ and $\varepsilon_0
> 0$ such that
\[
\mathrm{P}(\rho_T(\widetilde{X}^{x,\varepsilon}_t, \varphi) <
\delta) \geq \exp(-\varepsilon^{-2}[\overline{S}_{0,T}(\varphi) +
\gamma])
\]
for $\varepsilon < \varepsilon_0$ and $T > 0$, $\varphi \in
C([0,T], \mathbb{R}^d)$ such that $\varphi(0) = x$ and $T +
\overline{S}_{0,T}(\varphi) < C$.
\end{lemma}
\begin{lemma} \label{ler2} Suppose that
$\widetilde{a}^\varepsilon$, $\overline{a}$, ${\widetilde{b}}^\varepsilon$, and $\overline{b}$ are as above, and positive
constants  $k$, $K$, and $L$ are fixed. For $x \in \mathbb{R}^d$,
$T >0$, and $s \geq  0$, put
\[
\Phi(s) = \{ \varphi \in C([0,T], \mathbb{R}^d), \varphi(0) = x,
\overline{S}_{0,T}(\varphi) \leq s \}.
\]
For every $T > 0$, $\delta > 0$, $\gamma > 0$, and $s_0 > 0$, there
exist $\varkappa > 0$ and $\varepsilon_0 > 0$ such that for $x \in
\mathbb{R}^d$, $0 < \varepsilon \leq \varepsilon_0$ and $s \leq
s_0$, we have
\[
\mathrm{P}(\rho_T(\widetilde{X}^{x,\varepsilon}_t, \Phi(s)) \geq
\delta) \leq \exp(-\varepsilon^{-2}[s - \gamma]).
\]
\end{lemma}
Note that the choice of $\varkappa$ and $\varepsilon_0$ in the
above lemmas depends on the coefficients only through $k$, $K$ and $L$.

A sketch of the proof of these lemmas is provided in \cite{FKo}, so we'll not replicate it here. For the most part,
it is similar to the proof of the fact that
$\overline{S}_{0, T}(\varphi)$ serves as an action functional for the
process $\overline{X}^{x,\varepsilon}_t$ (see \cite{FWa}, \cite{DSt}).

We next state a corollary of the above two lemmas that will be
used in the paper. Define $D^\eta = \{ x \in D: {\rm dist}(x,
\partial D) \geq \eta \}$ and $T^\varepsilon(\lambda) = \exp(\lambda/\varepsilon^2)$.
Let
\[
v = \inf_{T, \varphi} \{ \overline{S}_{0,T}(\varphi): \varphi \in
C([0,T], \overline{D}), \varphi(0) = O, \varphi(T) \in \partial
D \}.
\]
Assume that
\[
(\widetilde{b}^\varepsilon(t, x), x-O),~ (\overline{b}(x), x-O)  \leq -c|x-O|^2
\]
whenever $x$ is in the $r$-neighborhood of $O$ and that
\[
 (\overline{b}(x), n(x)) \leq -c,~~~x \in \partial D.
\]
Moreover, we assume that
every solution of the equation $x'(t) = \overline{b}(x(t))$ that starts in $\overline{D}$ enters the $r$-neighborhood of $O$ in time that is
shorter than $c^{-1}$.

For $x \in D$, let $\widetilde{\tau}^{x,\varepsilon}$ be the first time when $\widetilde{X}^{x,\varepsilon}_t$ reaches the boundary of $D$. Thus$\widetilde{X}^{x,\varepsilon}_{\widetilde{\tau}^{x,\varepsilon} }$ is the location of the first
exit. Let
\[
\mathcal{A} = \{ x \in \partial D:
v = \inf_{T, \varphi} \{ \overline{S}_{0,T}(\varphi): \varphi \in
C([0,T], \overline{D}), \varphi(0) = O, \varphi(T) =x \} \}.
\]

\begin{lemma} \label{mlml}
Suppose that
$\widetilde{a}^\varepsilon$, $\overline{a}$, ${\widetilde{b}}^\varepsilon$, and $\overline{b}$ are as above. For each $\delta, \eta > 0$ there
are $\varkappa > 0$ and a function $\rho(\varepsilon)$ (that depend on $\widetilde{a}^\varepsilon$, $\overline{a}$, ${\widetilde{b}}^\varepsilon$, and $\overline{b}$ through $k$,
$K$, $L$, $c$, and $r$) such that $\lim_{\varepsilon \downarrow 0} \rho(\varepsilon)
= 0$ and
\\

(A)~~~~~~~~~~~~~~~~~~~~~~ $\mathrm{P}( \widetilde{\tau}^{x,\varepsilon}
 \leq T^\varepsilon(v + \delta))  \geq 1 -
\rho(\varepsilon)$ for $x \in D$,

(B)~~~~~~~~~~~~~~~~~~~~~ $ \mathrm{P}(\widetilde{\tau}^{x,\varepsilon} \geq
T^\varepsilon(v - \delta) ) \geq 1 - \rho(\varepsilon)$ for $x \in D^\eta$,

(C)~~~~~~~~~~~~~~~~~~~~~~~~~~~$ \mathrm{P}({\rm
dist}(\widetilde{X}^{x,\varepsilon}_{\widetilde{\tau}^{x,\varepsilon} }, \mathcal{A}) \leq
\eta) \geq 1 - \rho(\varepsilon)$ for $x \in D^\eta$.

(D)~~~~~~~~~~~~~~~~$\mathrm{P}(|\widetilde{X}^{x,\varepsilon}_t -O| < \eta) \geq 1 - \rho(\varepsilon)$ for $x \in D^\eta, t \in [T^\varepsilon(\delta), T^\varepsilon(v-\delta)]$,
\\
\\
provided that
\begin{equation} \label{clnss}
 \sup_{(t,x) \in \mathbb{R}^+
\times D} | \widetilde{b}_{i}^\varepsilon(t, x)-
\overline{b}_{i}(x)| \leq \varkappa,~~~~ \sup_{(t,x) \in \mathbb{R}^+
\times D }| \widetilde{a}_{ij}^\varepsilon(t, x)-
{\overline{a}}_{ij}(x)| \leq \varkappa.
\end{equation}
\end{lemma}
This lemma can be easily proved using a modification of Theorems
4.2 and 4.3 from Chapter 4 of \cite{FW} if we substitute our
Lemmas \ref{ler1} and \ref{ler2} for the corresponding results
concerning the case of time-independent coefficients.

An easy corollary of this lemma is that at an exponential time the
process either can be found in a small neighborhood of $O$ or
has earlier crossed the boundary of the domain.
\begin{lemma} \label{co1}
Suppose that
$\widetilde{a}^\varepsilon$, $\overline{a}$, ${\widetilde{b}}^\varepsilon$, and $\overline{b}$ are as above. For each $\delta, \eta > 0$ such that $O \in D^\eta$, there
is a function $\rho(\varepsilon)$ (that depends on $\widetilde{a}^\varepsilon$, $\overline{a}$, ${\widetilde{b}}^\varepsilon$, and $\overline{b}$ through $k$,
$K$, $L$, $c$, and $r$) such that $\lim_{\varepsilon \downarrow 0} \rho(\varepsilon)
= 0$ and
\[
\mathrm{P}(\widetilde{X}^{x,\varepsilon}_t \in D^\eta~~{\it
or}~~\widetilde{\tau}^{x,\varepsilon} \leq t) \geq 1 - \rho(\varepsilon)~~{\it for}~x \in D, t \in [T^\varepsilon(\delta), \infty).
\]
\end{lemma}
\proof Let $\eta_1 > 0$ be sufficiently small so that there is a
domain $\widetilde{D}$ with smooth boundary  such that
$\widetilde{D}^{\eta_1} = D$. Note that if $\eta_1$ is sufficiently small, then  Lemma~\ref{mlml} is applicable
to the domain $\widetilde{D}$. If the process does not reach
$\partial D$ by the time $t - T^\varepsilon(\delta)$, then we
can apply part (D) of Lemma~\ref{mlml} (with sufficiently small $\delta$) to the  domain
$\widetilde{D}$ and the process starting at
$X^{x,\varepsilon}_{t-T^\varepsilon(\delta)}$, and the result
follows from the Markov property. \qed

\section{Proof of the main result} \label{prml}
\begin{figure}     
\centerline{\includegraphics[height=3.2in, width= 4.3in,angle=0]{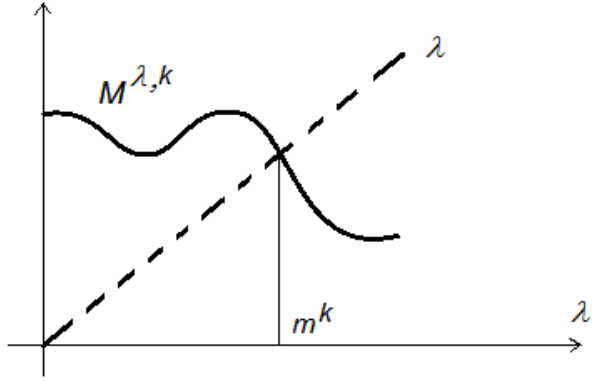}}
  \vskip -60pt
  \caption{Graph of $M^{\lambda,k}$.}
  \label{grm}
\end{figure}

In this section we prove Lemma~\ref{mnne}, which, as we discussed above, implies Theorem~\ref{mtee}. Let us first assume that $z_\lambda \equiv k$ (the general case will be considered below). An example of the graph of $M^{\lambda, k}$ is shown in Figure \ref{grm}.


Let $u^{\varepsilon,k}(t,x) = \mathrm{P}(\tau^{x, \varepsilon,k} \leq t)$.
Take an arbitrary $\eta > 0$ such that $O \in D^\eta$.
We claim that $ \lim_{\varepsilon \downarrow 0}  u^{\varepsilon,k}(T^\varepsilon(\lambda),x) = 0$ for each $\lambda < m^k$ uniformly in $x \in D^\eta$. Assume the contrary. Choose $\delta > 0$ and  $0 \leq \lambda_1 < \lambda_2 < m^k$ such that
\begin{equation} \label{imm}
\limsup_{\varepsilon \downarrow 0} \sup_{x \in D^\eta} u^{\varepsilon,k}(T^\varepsilon(\lambda_1),x) < \limsup_{\varepsilon \downarrow 0} \sup_{x \in D^\eta} u^{\varepsilon,k}(T^\varepsilon(\lambda_2),x),
\end{equation}
while
\begin{equation} \label{dlll}
\lambda_2 < M^{\lambda_1, k} - \delta.
\end{equation}
Let
\begin{equation} \label{coe1}
\overline{a}(x) = a(x, \lambda_1, k),~~\overline{b}(x) = b(x, \lambda_1, k),
\end{equation}
\begin{equation} \label{coe2}
\widetilde{a}(t,x) = a(x, \varepsilon^2 \ln( t) - T^\varepsilon(\lambda_1), k),~~\widetilde{b}(t,x) = b(x, \varepsilon^2 \ln( t)- T^\varepsilon(\lambda_1), k),
\end{equation}
$\widetilde{X}^{x,\varepsilon}_t$ be the corresponding process, and $\widetilde{\tau}^{x,\varepsilon}$ be the corresponding stopping time.
Choose $\varkappa > 0$ such that the conclusion of Lemma~\ref{mlml} holds. Due to the continuity assumptions
on $a$ and $b$, estimate (\ref{clnss}) holds if we restrict $t$ to the interval $[0, T^\varepsilon(\lambda_2) - T^\varepsilon(\lambda_1)]$ and choose $\lambda_2 - \lambda_1$ sufficiently small. Note that we can achieve (\ref{imm}), (\ref{dlll}) even with the requirement that $\lambda_2-\lambda_1$ is smaller than
any prescribed positive number.

By the definition of $u^{\varepsilon,k}$ and the Markov property of the process,
\begin{equation} \label{yyz}
u^{\varepsilon,k}(T^\varepsilon(\lambda_2),x) =
u^{\varepsilon,k}(T^\varepsilon(\lambda_1),x) + \mathrm{P}(\tau^{x, \varepsilon,k} > T^\varepsilon(\lambda_1), \widetilde{\tau}^{X^{x,\varepsilon,k}_{ T^\varepsilon(\lambda_1)}, \varepsilon} \leq T^\varepsilon(\lambda_2) -T^\varepsilon(\lambda_1))
\end{equation}
Observe that $T^\varepsilon(\lambda_2) - T^\varepsilon(\lambda_1) < T^\varepsilon(M^{\lambda_1, k}-\delta)$. Therefore, by part (B) of Lemma~\ref{mlml} and Lemma~\ref{co1}, the second term in the right hand side of (\ref{yyz}) tends to zero as $\varepsilon \downarrow 0$ uniformly in $x \in D^\eta$.  This, however, is a contradiction with (\ref{imm}), so we do have $ \lim_{\varepsilon \downarrow 0} u^{\varepsilon,k}(T^\varepsilon(\lambda),x) = 0$ for each $\lambda < m^k$ uniformly in $x \in D^\eta$.

Next, we claim that $ \lim_{\varepsilon \downarrow 0} u^{\varepsilon,k}(T^\varepsilon(\lambda),x) = 1$ for each $\lambda > m^k$ uniformly in $x \in D^\eta$. Since $m^k$ is the unique solution of the equation $M^{\lambda, k} = \lambda$, while $M^{\lambda, k}$ is bounded (as the coefficients are bounded), we have $M^{\lambda, k} < \lambda$ for $\lambda > m^k$. Choose $\delta > 0$ and $m^k < \lambda_1< \lambda_2$ such that
\[
M^{\lambda_1,k} + \delta < \lambda_1.
\]
We can achieve this while choosing $\lambda_2$ arbitrarily close to $m^k$. With
$\overline{a}$, $\overline{b}$,
$\widetilde{a}$, $\widetilde{b}$ defined by (\ref{coe1}) and (\ref{coe2}), choose $\varkappa > 0$ such that the conclusion of Lemma~\ref{mlml} holds. Due to the continuity assumptions
on $a$ and $b$, estimate (\ref{clnss}) holds if we restrict $t$ to the interval $[0, T^\varepsilon(\lambda_2) - T^\varepsilon(\lambda_1)]$ and choose $\lambda_2$ sufficiently close to $\lambda_1$. With these new $\lambda_1$ and $\lambda_2$, we again employ (\ref{yyz}), writing it as
\[
u^{\varepsilon,k}(T^\varepsilon(\lambda_2),x) =
\]
\[
u^{\varepsilon,k}(T^\varepsilon(\lambda_1),x) +
\mathrm{P}(\widetilde{\tau}^{X^{x,\varepsilon,k}_{ T^\varepsilon(\lambda_1)}, \varepsilon} \leq T^\varepsilon(\lambda_2) -T^\varepsilon(\lambda_1) |\tau^{x, \varepsilon,k} > T^\varepsilon(\lambda_1)) ( 1 - u^{\varepsilon,k}(T^\varepsilon(\lambda_1),x) ).
\]
The conditional probability in the right hand side tends to one uniformly in $x \in D^\eta$ as
$\varepsilon \downarrow 0$ by part (A) of Lemma~\ref{mlml} since $ T^\varepsilon(\lambda_2) - T^\varepsilon(\lambda_1) >  T^\varepsilon(\lambda_1) > T^\varepsilon( M^{\lambda_1,k} + \delta)$ for all sufficiently small $\varepsilon$. Therefore, $ \lim_{\varepsilon \downarrow 0} \inf_{x \in D^\eta} u^{\varepsilon,k}(T^\varepsilon(\lambda_2),x) = 1$. Since $\lambda_2$ can be chosen to be arbitrarily close to $m^k$ and since $u^{\varepsilon,k}(t,x)$ increases with $t$, we conclude that  $ \lim_{\varepsilon \downarrow 0} u^{\varepsilon,k}(T^\varepsilon(\lambda),x) = 1$ uniformly in $x \in D^\eta$.

Now let us prove that the distribution of $X^{x,\varepsilon,k}_{\tau^{x, \varepsilon,k}}$ converges to $x^k$. From our assumptions on the continuity of the coefficients it follows that for each $\lambda$ that is sufficiently close to $m^k$ the infimum $ \inf_{x \in \partial D}  V^{\lambda,k}(x)$ is achieved in a single point of the boundary that will be denoted by $x^k(\lambda)$. Moreover, $x^k(\lambda)$ is continuous at $m^k$. Given $\eta > 0$, choose $\lambda_1 < m^k$ such that ${\rm dist}(x^k(\lambda), x^k) < \eta$ for $\lambda \in [\lambda_1, m^k]$. With
$\overline{a}$, $\overline{b}$,
$\widetilde{a}$, $\widetilde{b}$ defined by (\ref{coe1}) and (\ref{coe2}), choose $\lambda_1 < m^k$ and $\lambda_2 > m_k$ in such a way that
(\ref{clnss}) holds if we restrict $t$ to the interval $[0, T^\varepsilon(\lambda_2) - T^\varepsilon(\lambda_1)]$, where $\varkappa$ is sufficiently small for the conclusion of Lemma~\ref{mlml} to hold (this may require modifying the previously selected $\lambda_1$ by making it larger).

As we showed above, $X^{x,\varepsilon,k}_t$ does not exit $D$ prior to time $T^\varepsilon(\lambda_1)$ with probability that tends to one uniformly in $x \in D^\eta$. By Lemma~\ref{co1},  $X^{x,\varepsilon,k}_{T^\varepsilon(\lambda_1)}$ belongs to $D^\eta$ with probability that tends to one. Therefore, by the above construction and part (C) of Lemma~\ref{mlml}, the process $X^{x,\varepsilon,k}_t$ reaches the boundary of $D$ for the first time in an $\eta$-neighborhood of $x^k(\lambda_1)$  with probability that tends to one. Since $\eta$ was arbitrary, this proves that the distribution of $X^{x,\varepsilon,k}_{\tau^{x, \varepsilon,k}}$ converges to $x^k$.

It remains now to get rid of the assumption that $z_\lambda \equiv k$.  Assume now that $z_\lambda$ has finitely many jumps at points $\lambda_1,...,\lambda_k$ and that none of the jumps happens at any of the points $m_1,...,m_s$. Suppose that $\sigma^z \geq \lambda_1$. As we showed above, $\mathrm{P}(\tau^{x, \varepsilon, z} > T^\varepsilon(\lambda_1)) \rightarrow 1$ as $\varepsilon \downarrow 0$. By Lemma~\ref{co1},  $X^{x,\varepsilon,z}_{T^\varepsilon(\lambda_1)}$ belongs to  $D^\eta$, with sufficiently small $\eta$,  with probability that tends to one. Using the Markov property of the process, we can now replace the values of $z$ on the interval $[0, \lambda_1)$ by $z_{\lambda_2}$, i.e., $z$ now has one fewer jump. At the same time, the coefficients of the process need to be set to be equal to zero on the interval $[0, T^\varepsilon(\lambda_1)]$. The earlier arguments still apply in this situation. Continuing by induction on the number of the jump, we obtain the desired result.

\section{The case of no repairing mechanism near the boundary} \label{exa}

Let us modify some of the assumptions on the coefficients $b$ and $\sigma$. Namely, instead of assuming that $b$ is Lipschitz continuous, we assume that it can be continued from $ \overline{D}$ to the entire space as a Lipschitz continuous function, yet
is equal to zero outside $ \overline{D}$. Inside $\overline{D}$, the drift is assumed to have the same properties as before. For simplicity, let us assume that $\sigma$ is an identity matrix and that $b$ does not depend on $\xi$. Thus the process under consideration is
\[
d X^{x,\varepsilon}_t = b(  X^{x,\varepsilon}_t, \varepsilon^2 \ln t) d t+\varepsilon
d W_t,~~X^{x,\varepsilon}_1 = x \in
\mathbb{R}^d.
\]
Let $D_1$ be a domain with smooth boundary such that $\overline{D} \subset D_1$. Let $\tau^{x, \varepsilon}_1$ be the first time when $X^{x, \varepsilon}_t$ reaches $\partial D_1$. We'll be interested in the  limiting distribution of $X^{x, \varepsilon}_{\tau^{x, \varepsilon}_1}$. The fact that $b = 0$ in $D_1 \setminus \overline{D}$ can be viewed as the absence of a repair mechanism in a part of the domain.

As follows from Theorem~\ref{mtee}, there is $x^* \in \partial D$ such that $\lim_{\varepsilon \downarrow 0} X^{x,\varepsilon}_{\tau^{x, \varepsilon}} = x^*$ in distribution for each $x \in D$. Note that $x^*$ is now non-random. In \cite{FK2} we considered the behavior of the process in $ \mathbb{R}^d \setminus D$ obtained from $X^{x,\varepsilon}_t$ by running the clock only when $X^{x,\varepsilon}_t$  is outside~${D}$. (The ageing mechanism was
not considered there, which is not a problem since the process in $ \mathbb{R}^d \setminus D$ is relevant now only at one time scale associated with the exit of  $X^{x,\varepsilon}_t$ from $D$). An adaptation of the results from \cite{FK2} gives us the following fact: for $\gamma \subset \partial D_1$ and $x \in D_1 \setminus D$,
\[
\lim_{\varepsilon \downarrow 0} \mathrm{P}(X^{x, \varepsilon}_{\tau^{x, \varepsilon}_1} \in \gamma ) = u^\gamma(x),
\]
where $u^\gamma$ is the unique solution of the following non-standard boundary problem
\[
\Delta u^\gamma (x) = 0,~~x\in D_1 \setminus \overline{D},
\]
\[
u^\gamma(x) = \begin{cases} 1, & x \in \gamma, \\ 0, & x \in \partial D_1 \setminus \gamma, \end{cases}
\]
\[
u^\gamma(x)~~{\rm is}~ {\rm constant}~{\rm on}~\partial D,
\]
\[
\langle \nabla u^\gamma(x^*), n(x^*) \rangle = 0,
\]
where $n(x)$ is the normal to $\partial D$ at $x$. For $x \in {D}$, we have $\lim_{\varepsilon \downarrow 0} \mathrm{P}(X^{x, \varepsilon}_{\tau^{x, \varepsilon}_1} \in \gamma ) = u^\gamma(x^*)$.
\\
\\


\noindent {\bf \large Acknowledgements}: While working on this
article, M. Freidlin was supported by NSF grant DMS-1411866
and L. Koralov was supported by NSF grant DMS-1309084.
\\
\\

\end{document}